\documentclass[10pt,leqno,twoside]{amsart}

\usepackage{amsmath,amssymb,amsfonts,amsthm}
\usepackage{fancyhdr}
\usepackage{amssymb,amsbsy}

\newtheorem{theorem}{\indent Theorem}[section]

\newtheorem{lemma}{\indent Lemma}[section]

\newcommand\be{\mathbb E}
\newcommand\bp{\mathbb P}

\title[Lower levels of the uniform
  recursive tree]{Degree distribution in the lower levels of the uniform
  recursive tree}

\author{\'Agnes Backhausz}
\address{Department of Probability Theory and Statistics\\Faculty of Science\\
 E\"otv\"os Lor\'and University\\P\'azm\'any P.~s.\ 1/C, H-1117
 Budapest, Hungary} 
\email{agnes@cs.elte.hu}
\author{Tam\'as F.~M\'ori}
\address{Department of Probability Theory and Statistics\\Faculty of Science\\
E\"otv\"os Lor\'and University\\P\'azm\'any P.~s. 1/C, 
H-1117 Budapest, Hungary}
\email{moritamas@ludens.elte.hu}
\dedicatory{\upshape Department of Probability Theory and Statistics,
Faculty of Science\\ E\"otv\"os Lor\'and University\\ P\'azm\'any P.~s. 1/C,
H-1117 Budapest, Hungary\\                         
\textit{E-mail address:} \texttt{agnes@cs.elte.hu, moritamas@ludens.elte.hu} \\
}


\thanks{The European Union and the European Social Fund have
provided financial support to the project under the grant agreement no.\ 
T\'AMOP 4.2.1./B-09/KMR-2010-0003.}

\begin{document}
\maketitle

\begin {abstract}
In this note we consider the $k$th level of the uniform random
recursive tree after $n$ steps, and prove that the proportion of nodes
with degree greater than $t\log n$ converges to $(1-t)^k$ almost
surely, as $n\to\infty$, for every $t\in(0,1)$. In addition, we show
that the number of degree $d$ nodes in the first level is
asymptotically Poisson distributed with mean $1$; moreover, they are
asymptotically independent for $d=1,2,\dots$.
\end{abstract}

\section{Introduction}

Let us consider the following random graph model. We start from a
single node labelled with $0$. At the $n$th step we 
choose a vertex at random, with equal probability, and independently
of the past. Then a new node, vertex $n$, is added to the graph, and it is
connected to the chosen vertex. In this way a random tree, the so
called uniform recursive tree, is built.

This model has a long and rich history. Apparently, the first
publication where the uniform recursive tree appeared was
\cite{Tapia}. Since then a huge number of papers have explored the
properties of this simple combinatorial structure.
 
Recursive trees serve as probabilistic models for
system generation, spread of contamination of organisms, pyramid
scheme, stemma construction of philology, Internet interface map, 
stochastic growth of networks, and many other areas of application,
see \cite{FHN} for references. For a survey of probabilistic
properties of uniform recursive trees see \cite{Dr} or \cite{SM}. 
Among others, it is known that this random tree has an asymptotic
degree distribution, namely, the 
proportion of nodes with degree $d$ converges, as $n\to\infty$, to
$2^{-d}$ almost surely. Another important quantity is the maximal
degree, which is known to be  asymptotically equal to $\log_2 n$ \cite{DL}.
Considering our graph a rooted tree, we can
define the levels of the tree in the usual way: level $k$ is the set
$L_n(k)$ of the vertices that are of distance $k$ from vertex $0$, the
root. It is not hard to find the a.s.\ asymptotics of the size of level $k$
after step $n$; it is
\[
|L_n(k)|\sim \be|L_n(k)|\sim\frac{(\log n)^k}{k!},\quad k=1,2,\dots\,.
\]

Recursive trees on nodes $0,\,1,\,\dots,\,n-1$ can be transformed into
permutations $\sigma=(\sigma_1,\sigma_2,\dots,\sigma_n)$ in the
following recursive way. Start from the identity permutation
$\sigma=(1,2,\dots,n)$. Then, taking the nodes $1,\,2,\,\dots,\,n-1$
one after another, update the permutation by swapping $\sigma_{i+1}$
and $\sigma_{i+1-j}$ if node $i$ was connected to node $j<i$ at the
time it was added to the tree. It is easy to see that in this way a 
one-to-one correspondence is set between trees and permutations, and
the uniform recursive tree is transformed into a uniform random
permutation. 

Another popular recursive tree model is the so called plane oriented
recursive tree. It was originally proposed by Szyma\'nski \cite{Szym}, but it
got in the focus of research after the seminal paper of
Barab\'asi and Albert \cite{BA}. A non-oriented
version of it starts from a single edge, and at each step a new vertex
is added to the graph. The new vertex is then connected to one of the
old nodes at random; the other endpoint of the new edge is chosen from
the existing vertices with probability proportional to the
instanteneous degree of the node (preferential attachment). This can
also be done in such a way that we select an 
edge at random with equal probability, then choose one of its endpoints.
In this tree the proportion of degree $d$ nodes converges
to $\frac{4}{d(d+1)(d+2)}$ with probability $1$.

Katona has shown \cite{Kat} that the same degree distribution 
can be observed if one is confined to any of the largest levels. On
the other hand, if we only consider a fixed level, the asymptotic degree
distribution still exists, but it becomes different \cite{M}. This
phenomenon has been observed in other random graphs, too. A general result
of that kind has been published recently \cite{BM}.

In the present note we will investigate the lower levels of the 
uniform recursive tree. We will show that, unlike in many scale free
recursive tree models, no asymptotic degree distribution
emerges. Instead, for almost all nodes in the lower levels the degree
sequence grows to infinity at the same rate as the overall maximum of
degrees does. We also investigate the number of degree $d$ vertices in
the first level for $d=1,\,2,\,\dots$, and show that they are
asymptotically i.i.d.\ Poisson with mean $1$. 

\section{Nodes of high degree in the lower levels}

Let $\deg_n(i)$ denote the degree of node $i$ after step $n$ $(i\le
n)$. Further, let $Z_{n,k}(t)$ denote  the proportion of nodes in
level $k$ with degree greater than $t\log n$. Formally,
\[
Z_{n,k}(t)=\frac{1}{|L_n(k)|}\,\bigl|\left\{i\le n: i\in L_n(k),\
    \deg_n(i)> t\log n\right\}\bigr|.
\]
The main result of this section is the following theorem.
\begin{theorem}\label{T2.1}
For $k=1,2,\dots$ and $0<t<1$
\[
\lim_{n\to\infty}Z_{n,k}(t)=(1-t)^k\ \ a.s.
\]
\end{theorem}
For the proof we need some auxiliary lemmas, interesting in their own
right. 

Let the number $n$ of steps be fixed, and $1<i<n$. Firstly, we are
interested in $X=\deg_n(i)-1$. 
\begin{lemma}\label{L2.1}
Let $0<\varepsilon<t<1$. Then for every $i>n^{1-t+\varepsilon}$ we have
\[
\bp(X>t\log n)\le \exp\left(-\frac{\varepsilon^2}{2t}\,\log n\right).
\]
\end{lemma}
\proof
$X=I_{i+1}+I_{i+2}+\dots+I_n$, where $I_j=1$, if vertex $i$ gets
a new edge at step $n$, and $0$ otherwise. These indicators are
clearly independent and $EI_j=1/j$, hence 
\[
\be X=\frac{1}{i+1}+\dots+\frac{1}{n}\,.
\]
Let us abbreviate it by $s$. Clearly, 
\[
\log\frac{n}{i+1}\le s\le \log\frac{n}{i}\,.
\]
Let $a>s$, then
by \cite[Theorem A.1.12]{AS} we have
\[
\bp(X\ge a)\le\left(e^{\beta-1}\beta^{-\beta}\right)^s,
\]
where $\beta=a/s$. Hence
\begin{multline*}
\bp(X\ge a)\le e^{a-s}\Bigl(\frac{s}{a}\Bigr)^a
=e^{a-s}\Bigl(1-\frac{a-s}{a}\Bigr)^a\\
=\exp\left(a-s-a\Bigl(\frac{a-s}{a}+\frac 12\Bigl(\frac{a-s}{a}
\Bigr)^2+\dots \Bigr)\right)\le\exp\left(-\frac{(a-s)^2}{2a}\right).
\end{multline*}
Now, set $a=t\log n$. Then $s\le(t-\varepsilon)\log n$, and
\[
\bp(X\ge t\log n)\le\exp\left(-\frac{(t\log n-s)^2}{2t\log n}\right)
\le\exp\left(-\frac{\varepsilon^2}{2t}\,\log n\right).
\]
\qed

\begin{lemma}\label{L2.2}
Let $0<t<1$, and $0<\varepsilon<1-t$. Then for every
$i\le n^{1-t-\varepsilon}-1$ we have 
\[
\bp(X\le t\log n)\le \exp\left(-\frac{\varepsilon^2}{2(t+\varepsilon)}\,
\log n\right).
\]
\end{lemma}
\proof
This time $s>\log\frac{n}{i+1}\ge(t+\varepsilon)\log n$, thus
\cite[Theorem A.1.13]{AS} implies that
\[
\bp(X\le t\log n)\le\exp\left(-\frac{(s-t\log n)^2}{2s}\right).
\]
Notice that the exponent in the right-hand side, as a function
of s, is decreasing for $s>t\log n$. Therefore $s$ can be replaced
by $(t+\varepsilon)\log n$, and the proof is complete.
\qed

\indent \textbf{Proof of Theorem \ref{T2.1}}.
Since $\deg_n(i)$ is approximately equal to $\log\frac{n}{i}$, it
follows that $\deg_n(i)\ge t\log n$ is approximately equivalent to
$i\le n^{1-t}$. Basing on Lemmas \ref{L2.1} and \ref{L2.2} we can
quantify this heuristic reasoning. 

Let $0<\varepsilon<\min\{t,\,1-t\}$, and
$a=a(n)=\left\lfloor n^{1-t-\varepsilon}\right\rfloor-1$, 
$b=b(n)=\left\lceil n^{1-t+\varepsilon}\right\rceil$. Then
by Lemma \ref{L2.2}
\begin{align*}
&\bp\bigl(\exists i\in L_n(k)\text{ such that }i\le a,\ 
\deg_n(i)\le 1+t\log n\bigr)\\
&\qquad\le
\sum_{i=1}^a \bp\bigl(i\in L_n(k),\ \deg_n(i)\le 1+t\log n\bigr)\\
&\qquad =\sum_{i=1}^a \bp\bigl(i\in L_n(k)\bigr)\bp\bigl(\deg_n(i)\le
1+t\log n\bigr)\\ 
&\qquad\le \be L_n(k)\cdot\exp\left(-\frac{\varepsilon^2}{2(t+\varepsilon)}\,
\log n\right).
\end{align*}
Similarly, by Lemma \ref{L2.1},
\begin{align*}
&\bp\bigl(\exists i\in L_n(k)\text{ such that }i>b,\ 
\deg_n(i)> 1+t\log n\bigr)\\
&\qquad\le
\sum_{i=b+1}^n \bp\bigl(i\in L_n(k),\ \deg_n(i)> 1+t\log n\bigr)\\
&\qquad =\sum_{i=b+1}^n \bp\bigl(i\in L_n(k)\bigr)\bp\bigl(\deg_n(i)>
1+t\log n\bigr)\\ 
&\qquad\le \be L_n(k)\cdot\exp\left(-\frac{\varepsilon^2}{2t}\,
\log n\right).
\end{align*}
Introduce the events 
\[
A(n)=\bigl\{L_a(k)\subset\{i\in L_n(k):\deg_n(i)>1+t\log n\}
\subset L_b(k)\bigr\},
\]
then the probability of their complements can be estimated as follows.
\[
\bp\left(\overline{A(n)}\right)\le 2\be|L_n(k)|\,
\exp\left(-\frac{\varepsilon^2}{2t}\,\log n\right).
\]
Note that $|L_a(k)|\sim(1-t-\varepsilon)^k|L_n(k)|$, and 
$|L_b(k)|\sim(1-t+\varepsilon)^k|L_n(k)|$, a.s. 

Let $c>2(t+\varepsilon)\varepsilon^{-2}$, then
$\sum_{n=1}^\infty \bp\left(\overline{A(n^c)}\right)<\infty$, hence by
the Borel--Cantelli lemma it follows almost surely that $A(n^c)$
occurs for every $n$ large enough. Consequently, 
\begin{multline*}
(1-t-\varepsilon)^k\left|L_{n^c}(k)\right|\bigl(1+o(1)\bigr)\\
\le \left|\left\{i\in L_{n^c}(k):\deg_{n^c}(i)>1+t\log(n^c)\right\}
\right|\\
\le(1-t+\varepsilon)^k\left|L_{n^c}(k)\right|\bigl(1+o(1)\bigr).
\end{multline*}
This implies 
\[
\liminf_{n\to\infty} Z_{n^c,k}(t)\ge(1-t-\varepsilon)^k\text{ and }
\limsup_{n\to\infty} Z_{n^c,k}(t)\le(1-t+\varepsilon)^k
\]
for every positive $\varepsilon$, hence Theorem \ref{T2.1} is proven
along the subsequence $(n^c)$. 

To the indices in between we can apply
the following esimation. For $n^c\le N \le(n+1)^c$ with sufficiently
large $n$ we have
\begin{align*}
Z_{N,k}(t)&\le\frac{1}{\left|L_{n^c}(k)\right|}\,\left|\left\{i\in 
L_{(n+1)^c}(k):\deg_{(n+1)^c}(i)\ge t\log(n^c)\right\}\right|\\
&=\frac{\left|L_{(n+1)^c}(k)\right|}{\left|L_{n^c}(k)\right|}\,   
Z_{(n+1)^c,k}\left(t\,\frac{\log n}{\log(n+1)}\right).
\end{align*}
Here the first term tends to $1$, while the second term's asymptotic
behaviour is just the same as that of $Z_{(n+1)^c,k}(t)$. Hence
$Z_{N,k}(t)\le\bigl(1+o(1)\bigr)(1-t)^k$.

Similarly,
\begin{align*}
Z_{N,k}(t)
&\ge\frac{\left|L_{n^c}(k)\right|}{\left|L_{(n+1)^c}(k)\right|}\,   
Z_{n^c,k}\left(t\,\frac{\log(n+1)}{\log n}\right)\\
&=\bigl(1+o(1)\bigr)Z_{n^c,k}(t)\\
&=\bigl(1+o(1)\bigr)(1-t)^k.
\end{align*}
This completes the proof.
\qed

\section{Nodes of small degree in the first level}
Looking at the picture Theorem \ref{T2.1} shows us on the degree
distribution one can naturally ask how many points of fixed degree
remain in the lower levels at all. In this respect the first level and the
other ones behave differently. It is easy to see that 
degree $1$ nodes in level $1$ correspond to the fixed points of the random
permutation described in the Introduction. Hence their number has a
Poisson limit distribution with parameter $1$ without any
normalization. More generally, let 
\[
X[n,d]=\bigl|\{i\in L_n(1): \deg_n(i)=d\}\bigr|;
\]
this is the number of nodes with degree $d$
in the first level after $n$ steps. 

The main result of this section is the following limit theorem.
\begin{theorem}\label{T3.1}
$X[n,1],\,X[n,2],\,\dots$ are asymptotically i.i.d.\ Poisson with mean
$1$, as $n\to\infty$.
\end{theorem}

\proof
We will apply the method of moments in the following form.

For any real number $a$ and nonnegative integer $k$ let us define
$(a)_0=1$, and $(a)_k=a(a-1)\cdots(a-k+1)$, $k=1,2,\dots$\,. In order
to verify the limiting joint distribution in Theorem \ref{T3.1} it
suffices to show that
\begin{equation}\label{3.1}
\lim_{n\to\infty}\be\left(\prod_{i=1}^d\bigl(X[n,i]\bigr)_{k_i}\right)=1
\end{equation}
holds for every $d=1,2,\dots$, and nonnegative integers
$k_1,\,\dots,\,k_d$. This can easily be seen from the following
expansion of the joint probability generating function of the random
variables $Y[n,1],\,\dots,\,Y[n,d]$.
\[
\be\left(\prod_{i=1}^d z_i^{Y[n,i]}\right)=
\sum_{k_1=0}^\infty\dots\sum_{k_d=0}^\infty\be
\left(\prod_{i=1}^d\bigl(X[n,i]\bigr)_{k_i}\right)\prod_{i=1}^d
\frac{(z_i-1)^{k_i}}{k_i!}.
\]

In the proof we shall rely on the following obvious identities.
\begin{gather}
(a+1)_k-(a)_k=k\,(a)_{k-1},\label{3.2}\\
a\bigl[(a-1)_k(b+1)_\ell-(a)_k(b)_\ell\bigr]=\ell\,(a)_{k+1}
(b)_{\ell-1}-k\,(a)_k(b)_\ell,\label{3.3}\\
\sum_{a=k}^n(a)_k=\frac{1}{k+1}\,(n+1)_{k+1}.\label{3.4}
\end{gather}

Let us start from the conditional expectation of the quantity under
consideration with respect to the sigma-field generated by the
past of the process.
\begin{equation}\label{3.5}
\be\Biggl(\prod_{i=1}^d\bigl(X[n+1,i]\bigr)_{k_i}\Biggm|\mathcal F_n
\Biggr)=\prod_{i=1}^d\bigl(X[n,i]\bigr)_{k_i}+\sum_{j=0}^d S_j,
\end{equation}
where in the rightmost sum $j$ equals $0,\,1,\,\dots,\,d$, according
that the new vertex at step $n+1$ is 
connected to the root $(j=0)$, or to a degree $j$ node in level
$1$. This happens with (conditional) probability 
$\dfrac 1n\,,\,\dfrac{X[n,1]}{n}\,,\ \dots,\,\dfrac{X[n,d]}{n}\,,$
respectively. That is, 
\begin{equation*}
S_0=\frac 1n\,\prod_{i=2}^d\bigl(X[n,i]\bigr)_{k_i}\left[
\bigl(X[n,1]+1\bigr)_{k_1}-\bigl(X[n,1]\bigr)_{k_1}\right],
\end{equation*}
and for $1\le j\le d-1$
\begin{multline*}
S_j=\frac{X[n,j]}{n}\prod_{i\ne\{j,j+1\}}\bigl(X[n,i]\bigr)_{k_i}
\times\\
\times\Bigl[\bigl(X[n,j]-1\bigr)_{k_j}\bigl(X[n,j+1]+1\bigr)_{k_{j+1}}
-\bigl(X[n,j]\bigr)_{k_j}\bigl(X[n,j+1]\bigr)_{k_{j+1}}\Bigr].
\end{multline*}
Finally,
\begin{equation*}
S_d=\frac{X[n,d]}{n}\ \prod_{i=1}^{d-1}\bigl(X[n,i]\bigr)_{k_i}
\left[\bigl(X[n,d]-1\bigr)_{k_d}-\bigl(X[n,d]\bigr)_{k_d}\right].
\end{equation*}

Let us apply \eqref{3.2} to $S_0$ with $k=k_1$, \eqref{3.3} to $S_j$
with $k=k_j$, $\ell=k_{j+1}$ $(1\le j\le d-1)$, and \eqref{3.3} to
$S_d$ with $k=k_d$, $\ell=0$, to obtain
\begin{gather}
S_0=\frac{k_1}{n}\ \prod_{i=2}^d\bigl(X[n,i]\bigr)_{k_i}
\,\bigl(X[n,1]\bigr)_{k_1-1}\,,\label{3.9}\\
\begin{aligned}\label{3.10}
S_j=\frac{k_{j+1}}{n}\ \prod_{i\ne\{j,j+1\}}\bigl(X[n,i]\bigr)_{k_i}
\,\bigl(X[n,j]\bigr)_{k_j+1}\,\bigl(X[n,j+1]\bigr)_{k_{j+1}-1}\\
-\frac{k_j}{n}\ \prod_{i=1}^d\bigl(X[n,i]\bigr)_{k_i}\,,
\end{aligned}\\
S_d=-\frac{k_d}{n}\ \prod_{i=1}^d\bigl(X[n,i]\bigr)_{k_i}\,.
\label{3.11}
\end{gather}
In \eqref{3.9}--\eqref{3.10} it can happen that some of the $k_j$ are
zero, and, though $(a)_{-1}$ has not been defined, it always gets a zero
multiplier, thus the expressions do have sense.
Let us plug \eqref{3.9}--\eqref{3.11} into \eqref{3.5}.
\begin{multline*}
\be\Biggl(\prod_{i=1}^d\bigl(X[n+1,i]\bigr)_{k_i}\Biggm|\mathcal F_n
\Biggr)\\
=\prod_{i=1}^d\bigl(X[n,i]\bigr)_{k_i}\Biggl(1-\frac 1n\,
\sum_{j=1}^d k_j\Biggr)
+\frac{k_1}{n}\,\prod_{i=2}^d\bigl(X[n,i]\bigr)_{k_i}
\,\bigl(X[n,1]\bigr)_{k_1-1}\\
+\sum_{j=1}^{d-1}\,\frac{k_{j+1}}{n}\,\prod_{i\ne\{j,j+1\}}
\bigl(X[n,i]\bigr)_{k_i}\,\bigl(X[n,j]\bigr)_{k_j+1}\,
\bigl(X[n,j+1]\bigr)_{k_{j+1}-1}.
\end{multline*}
Introducing
\[
E(n,k_1,\dots,k_d)=\be\left(\prod_{i=1}^d\bigl(X[n,i]\bigr)_{k_i}\right),
\quad K=k_1+\dots+k_d,
\]
we have the following recursion.
\begin{gather*}
E(n+1,k_1,\dots,k_d)=\left(1-\frac{K}{n}\right)E(n,k_1,\dots,k_d)+
\frac{k_1}{n}\,E(n,k_1-1,k_2,\dots,k_d)\\
+\sum_{j=1}^{d-1}\,\frac{k_{j+1}}{n}\,E(n,k_1,\dots,k_j+1,k_{j+1}-1,\dots,k_d),
\end{gather*}
or equivalently,
\begin{multline}\label{3.12}
(n)_K E(n+1,k_1,\dots,k_d)=(n-1)_K E(n,k_1,\dots,k_d)\\
+(n-1)_{K-1}\sum_{j=1}^d k_j E(k_1,\dots,k_{j-1}+1,k_j-1,\dots,k_d).
\end{multline}
Based on \eqref{3.12}, the proof can be completed by induction on the
exponent vectors $(k_1,\dots,k_n)$. We say that $\underline
k=(k_1,k_2,\dots,k_d)$ is majorized by
$\underline{\ell}=(\ell_1,\ell_2,\dots,\ell_d)$, if $k_d\le\ell_d$,
$k_{d-1}+k_d\le\ell_{d-1}+\ell_d$, \dots,
$k_1+\dots+k_{d}\le\ell_1+\dots+\ell_{d}$. This is a total 
order on $\mathbb N^d$. 

Now, \eqref{3.1} is clearly holds for $\underline k=(1,0,\dots,0)$,
since $\be X[n,1]=1$ for every $n=1,2,\dots\,$, which is obvious
considering the fixed points of a random permutation.

In every term of the sum on the right hand
side of \eqref{3.12} the argument of $E(\,\cdot\,)$ is majorized by
$\underline k=(k_1,\dots,k_d)$, hence the induction hypothesis can be
applied to them. We get that
\[
(n)_K E(n+1,\underline k)=(n-1)_K E(n,\underline k)
+(n-1)_{K-1}K\bigl(1+o(1)\bigr),
\] 
from which \eqref{3.4} gives that $(n-1)_K E(n,\underline k)
\sim(n-1)_K$, that is,
\[
\lim_{n\to\infty}E(n,k_1,\dots,k_d)=1,
\]
as needed. \qed

Turning to higher levels one finds the situation changed. Fixing a
degree $d$ we find, roughly speaking, that each node in level $k-1$
has a Poisson number of degree $d$ children in level $k$ (a freshly
added node is considered the child of the old node it is connected
to). Now, strong-law-of-large-numbers-type heuristics imply that
the number of nodes with degree $d$ in level $k\ge 2$ is approximately
equal to $|L_n(k-1)|$, that is, their proportion is 
\[
\approx\frac{|L_n(k-1)|}{|L_n(k)|}\sim\frac{1}{k\log n}\,.
\]

Another interesting problem worth of dealing with is the number of
nodes with unusually high degree. In every fixed level Theorem
\ref{T2.1} implies that the proportion of nodes with degree higher
than $\log n$ is asymptotically negligible, but they must exist, since
the maximal degree is approximately $\log_2 n=\log_2e\,\log n$. How
many of them are there? We are planning to return to this issue in a
separate paper.

\end{document}